\documentclass[10pt]{amsart}
\usepackage{amsfonts}
\usepackage{ifthen}
\usepackage{amsthm}
\usepackage{amsmath}
\usepackage{graphicx}
\usepackage{amscd,amssymb,amsthm}
\usepackage{epstopdf}

\newcounter{minutes}
\divide\time by 60
\newcounter{hours}
\multiply\time by 60 \addtocounter{minutes}{-\time}

\newcommand{\real}{\operatorname{Re}}
\newcommand{\logo}{\operatorname{Log}}

\newtheorem{Theorem}{Theorem}


\title[\'A. Baricz, H. Orhan, R. Sz\'asz/The radius of $\alpha$-convexity of Bessel functions]{The radius of $\alpha$-convexity of normalized Bessel functions of the first kind}

\author[]{\'Arp\'ad Baricz}
\address{Department of Economics, Babe\c{s}-Bolyai University, 400591 Cluj-Napoca, Romania}
\address{Institute of Applied Mathematics, \'Obuda University, 1034 Budapest, Hungary}
\email{bariczocsi@yahoo.com}

\author[]{Halit Orhan}
\address{Department of Mathematics, Ataturk University, 25240 Erzurum, Turkey}
\email{horhan@atauni.edu.tr}

\author[]{R\'obert Sz\'asz}
\address{Department of Mathematics and Informatics, Sapientia Hungarian University of Transylvania, 540485 T\^argu-Mure\c{s}, Romania}
\email{rszasz@ms.sapientia.ro}

\begin{document}

\def\thefootnote{}
\footnotetext{ \texttt{File:~\jobname .tex,
          printed: \number\year-\number\month-0\number\day,
          \thehours.\ifnum\theminutes<10{0}\fi\theminutes}
} \makeatletter\def\thefootnote{\@arabic\c@footnote}\makeatother

\maketitle

\begin{abstract}
The radii of $\alpha$-convexity are deduced for three different kind of normalized Bessel functions of the first kind and it is shown that these radii are
between the radii of starlikeness and convexity, when $\alpha\in[0,1],$ and they are decreasing with respect to the parameter $\alpha.$ The results presented in this paper unify some recent results on the radii of starlikeness and convexity for normalized Bessel functions of the first kind. The key tools in the proofs are some interlacing properties of the zeros of some Dini functions and the zeros of Bessel functions of the first kind.
\end{abstract}

\section{\bf Introduction and the main results}

Let $\mathbb{D}(0,{r})$ be the open disk $\left\{ {z\in \mathbb{C}:\left\vert
z\right\vert <r}\right\} ,$ where $r>0,$ and set $\mathbb{D}=\mathbb{D}(0,{1})$. By $\mathcal{A}$ we mean the class of analytic functions $f:\mathbb{D}(0,r)\to
\mathbb{C}$ which satisfy the usual normalization conditions $f(0)=f'(0)-1=0.$ Denote by $\mathcal{S}$ the class of functions belonging to $\mathcal{A}$ which are univalent in $\mathbb{D}(0,r)$ and let $\mathcal{S}^{\ast }(\beta )$ be the subclass of $\mathcal{S}$ consisting of functions which are starlike of order $\beta $ in $\mathbb{D}(0,r),$ where $0\leq \beta <1.$ The analytic characterization of this class of functions is
\begin{equation*}
\mathcal{S}^{\ast }(\beta )=\left\{ f\in \mathcal{S}\ :\ \real \left(\frac{zf'(z)}{f(z)}\right)>\beta\ \ \mathrm{for\ all}\ \ z\in \mathbb{D}(0,r) \right\},
\end{equation*}
while the real number
\begin{equation*}
r_{\beta }^{\ast}(f)=\sup \left\{ r>0\ :\ \real \left(\frac{zf'(z)}{f(z)}\right)>\beta\ \ \mathrm{\ for\ all}\ \ z\in \mathbb{D}(0,r)\right\}
\end{equation*}
is called the radius of starlikeness of order $\beta $ of the
function $f.$ Note that $r^{\ast }(f)=r_0^{\ast}(f)$ is the largest radius
such that the image region $f(\mathbb{D}(0,{r^{\ast }(f)}))$ is a starlike
domain with respect to the origin. Also, let $\mathcal{K}(\beta)$ be the subclass of $\mathcal{S}$ consisting
of functions which are convex of order $\beta $ in $\mathbb{D}(0,r),$ where $0\leq \beta <1.$ The well-known analytic characterization of this class of functions is
\begin{equation*}
\mathcal{K}(\beta )=\left\{ f\in \mathcal{S}\ :\ \real \left(1+\frac{zf''(z)}{f'(z)}\right)>\beta\ \ \mathrm{for\ all}\ \ z\in\mathbb{D}(0,r) \right\},
\end{equation*}
and the real number
\begin{equation*}
r_{\beta }^{c}(f)=\sup \left\{ r>0\ :\ \real \left(1+\frac{zf''(z)}{f'(z)}\right)>\beta\ \ \mathrm{\ for\ all}\ \ z\in \mathbb{D}(0,r)\right\}
\end{equation*}
is called the radius of convexity of order $\beta $ of the
function $f.$ Note that $r^{c}(f)=r_0^{c}(f)$ is the largest radius such
that the image region $f(\mathbb{D}(0,{r^{c}(f)}))$ is a convex domain with
respect to the origin. Furthermore, let $\mathcal{M}(\alpha,\beta)$ be the subclass of $\mathcal{S}$
consisting of functions which are $\alpha-$convex of order $\beta $ in $\mathbb{D}(0,r),$ where $\alpha\in\mathbb{R}$ and $0\leq \beta <1.$ The
analytic characterization of this class of functions is
\begin{equation*}
\mathcal{M}(\alpha,\beta )=\left\{ f\in \mathcal{S}\ :\ \real \left((1-\alpha)\frac{zf'(z)}{f(z)}+\alpha\left(1+\frac{zf''(z)}{f'(z)}\right)\right)>\beta\ \ \mathrm{for\ all}\ \ z\in \mathbb{D}(0,r)\right\},
\end{equation*}
while the real number
\begin{equation*}
r_{\alpha,\beta }(f)=\sup \left\{ r>0\ :\ \real \left((1-\alpha)\frac{zf'(z)}{f(z)}+\alpha\left(1+\frac{zf''(z)}{f'(z)}\right)\right)>\beta\ \ \mathrm{\ for\ all}\ \ z\in \mathbb{D}(0,r)\right\}
\end{equation*}
is called the radius of $\alpha-$convexity of order $\beta $ of the function $f.$ The radius of $\alpha-$convexity of order $\beta$ is the generalization of
the radius of starlikeness of order $\beta$ and of the radius of convexity of order $\beta.$ We have $r_{0,\beta}(f)=r_{\beta }^{\ast}(f)$ and $r_{1,\beta}(f)=
r_{\beta }^{c}(f).$ For more details on starlike, convex and $\alpha$-convex functions we refer to \cite{duren,mocanu,moc} and to the references therein.

The Bessel function of the first kind of order $\nu$ is defined by \cite[p. 217]{nist}
$$J_{\nu}(z)=\sum_{n\geq0}\frac{(-1)^{n}}{n!\Gamma(n+\nu+1)}\left(\frac{z}{2}\right)^{2n+\nu},\ \ z\in\mathbb{C}.$$
In this paper we focus on the following  normalized forms
$$f_{\nu}(z)=\left(2^{\nu}\Gamma(\nu+1)J_{\nu}(z)\right)^{\frac{1}{\nu}}=z-\frac{1}{4\nu(\nu+1)}z^3+\dots,\ \nu\neq0,$$
$$g_{\nu}(z)=2^{\nu}\Gamma(\nu+1)z^{1-{\nu}}J_{\nu}(z)=z-\frac{1}{4(\nu+1)}z^3+\frac{1}{32(\nu+1)(\nu+2)}z^5-\dots,$$
$$h_{\nu}(z)=2^{\nu}\Gamma(\nu+1)z^{1-\frac{\nu}{2}}J_{\nu}(\sqrt{z})=z-\frac{1}{4(\nu+1)}z^2+\dots,$$
where $\nu>-1.$ We note that in fact for $z\in\mathbb{C}\setminus\{0\}$ we have $$f_{\nu}(z)=\exp\left(\frac{1}{\nu}\logo\left(2^{\nu}\Gamma(\nu+1)J_{\nu}(z)\right)\right),$$
where $\logo$ represents the principal branch of the logarithm, and in this paper every
multi-valued function is taken with the principal branch. We also mention that the univalence, starlikeness and convexity of Bessel function of the first kind were studied extensively in several papers. We refer to \cite{mathematica,publ,lecture,bsk,samy,anal,brown,brown2,brown3,todd,szasz,szasz2} and to the references therein.

In this paper we make a further contribution to the subject by showing the following new sharp results contained in Theorems \ref{th1}, \ref{th2} and \ref{th3}. The proofs of Theorems \ref{th1}, \ref{th2} and \ref{th3} can be found in section 2. 

\begin{Theorem}\label{th1}
If $\nu>0,$ $\alpha\geq0$ and $\beta\in[0,1),$ then the radius of
$\alpha$-convexity of order $\beta$ of the function $f_\nu$ is the smallest positive
root of the equation
$$
\alpha\left(1+\frac{rJ_\nu''(r)}{J_\nu'(r)}\right)+\left(\frac{1}{\nu}-\alpha\right)\frac{rJ_\nu'(r)}{J_\nu(r)}=\beta.
$$
The radius of $\alpha$-convexity satisfies $r_{\alpha,\beta}(f_{\nu})\leq j_{\nu,1}'<j_{\nu,1},$ where $j_{\nu,1}$ and $j_{\nu,1}'$ denote the first positive zeros of $J_{\nu}$ and $J_{\nu}',$ respectively. Moreover, the function $\alpha\mapsto r_{\alpha,\beta}(f_{\nu})$ is strictly decreasing on $[0,\infty)$ and consequently we have $r_{\beta}^{c}(f_{\nu})<r_{\alpha,\beta}(f_{\nu})<r_{\beta }^{\ast}(f_{\nu})$ for all $\alpha\in(0,1),$ $\beta\in[0,1)$ and $\nu>0.$
\end{Theorem}

\begin{figure}[!ht]
   \centering
       \includegraphics[width=10cm]{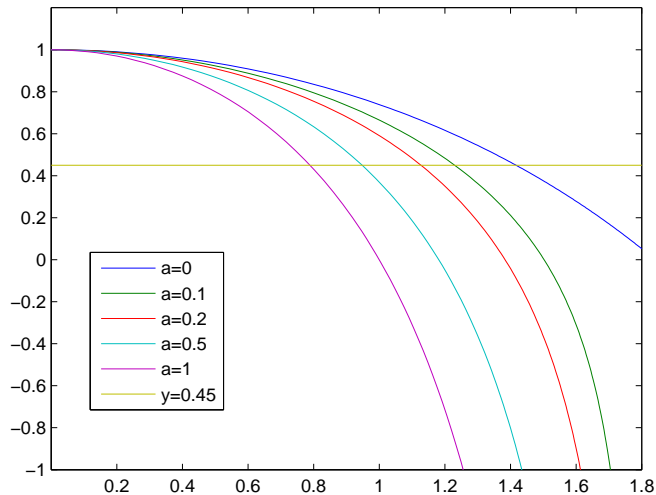}
       \caption{The graph of the function $r\mapsto \mathbb{J}(a,f_1(r))$ for $a\in\{0,0.1,0.2,0.5,1\}$ on $[0,1.8].$}
       \label{fig1}
\end{figure}

It is worth to mention that the cases $\alpha=0$ and $\alpha=1$ of the above Theorem were considered recently in \cite[Theorem 1(a)]{bsk} and \cite[Theorem 1.1]{anal}. Our Theorem \ref{th1} is a common generalization of these results. Figure \ref{fig1} illustrates the fact that if $\alpha\in[0,1],$ then the radius of $\alpha$-convexity of the function $f_{\nu}$ is between its radii of convexity and starlikeness, that is, $r_{\beta}^{c}(f_{\nu})<r_{\alpha,\beta}(f_{\nu})<r_{\beta }^{\ast}(f_{\nu})$ for all $\alpha\in(0,1),$ $\beta\in[0,1)$ and $\nu>0.$ We considered the particular cases when $\beta=0.45,$ $\nu=1$ and $\alpha\in\{0,0.1,0.2,0.5,1\}.$

\begin{Theorem}\label{th2}
If $\nu>-1,$ $\alpha\geq0$ and $\beta\in[0,1),$ then the radius of
$\alpha$-convexity of order $\beta$ of the function $g_{\nu}$ is the smallest positive root of
the equation
$$1+(\alpha-1)\frac{rJ_{\nu+1}(r)}{J_{\nu}(r)}+\alpha r\frac{rJ_{\nu+2}(r)-3J_{\nu+1}(r)}{J_{\nu}(r)-rJ_{\nu+1}(r)}=\beta.$$
The radius of $\alpha$-convexity satisfies $r_{\alpha,\beta}(g_{\nu})\leq\alpha_{\nu,1}<j_{\nu,1},$ where $\alpha_{\nu,1}$ is the first positive zero of the Dini function $z\mapsto (1-\nu)J_\nu(z)+zJ'_\nu(z).$ Moreover, the function $\alpha\mapsto r_{\alpha,\beta}(g_{\nu})$ is strictly decreasing on $[0,\infty)$ and consequently we have $r_{\beta}^{c}(g_{\nu})<r_{\alpha,\beta}(g_{\nu})<r_{\beta }^{\ast}(g_{\nu})$ for all $\alpha\in(0,1),$ $\beta\in[0,1)$ and $\nu>-1.$
 \end{Theorem}

\begin{figure}[!ht]
   \centering
       \includegraphics[width=10cm]{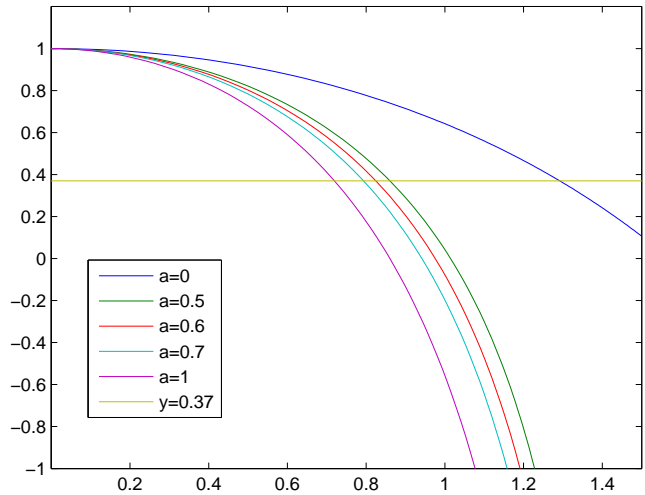}
       \caption{The graph of the function $r\mapsto \mathbb{J}(a,g_{0.5}r))$ for $a\in\{0,0.5,0.6,0.7,1\}$ on $[0,1.5].$}
       \label{fig2}
\end{figure}

We also note that the cases $\alpha=0$ and $\alpha=1$ of the above Theorem were considered recently in \cite[Theorem 1(b)]{bsk} and \cite[Theorem 1.2]{anal}. Our Theorem \ref{th2} is a common generalization of these results. Figure \ref{fig2} illustrates the fact that when $\alpha\in[0,1]$ the radius of $\alpha$-convexity of the function $g_{\nu}$ is between its radii of convexity and starlikeness, that is, $r_{\beta}^{c}(g_{\nu})<r_{\alpha,\beta}(g_{\nu})<r_{\beta }^{\ast}(g_{\nu})$ for all $\alpha\in(0,1),$ $\beta\in[0,1)$ and $\nu>-1.$ We considered the particular cases when $\beta=0.37,$ $\nu=0.5$ and $\alpha\in\{0,0.5,0.6,0.7,1\}.$

\begin{Theorem}\label{th3}
If $\nu>-1,$ $\alpha\geq0$ and $\beta\in[0,1),$ then the radius of
$\alpha$-convexity of order $\beta$ of the function $h_{\nu}$ is the smallest positive root of
the equation
$$(1-\alpha)\left(1-\frac{r^{\frac{1}{2}}}{2}\cdot\frac{J_{\nu+1}(r^{\frac{1}{2}})}{J_{\nu}(r^{\frac{1}{2}})}\right)+\alpha\left(1+\frac{r^{\frac{1}{2}}}{2}\cdot\frac{r^{\frac{1}{2}}J_{\nu+2}(r^{\frac{1}{2}})-
 4J_{\nu+1}(r^{\frac{1}{2}})}{2J_{\nu}(r^{\frac{1}{2}})-r^{\frac{1}{2}}J_{\nu+1}(r^{\frac{1}{2}})}\right)=\beta.$$
The radius of $\alpha$-convexity satisfies $r_{\alpha,\beta}(h_{\nu})\leq\beta_{\nu,1}^2<j_{\nu,1}^2,$ where $\alpha_{\beta,1}$ is the first positive zero of the Dini function $z\mapsto (2-\nu)J_\nu(z)+zJ'_\nu(z).$ Moreover, the function $\alpha\mapsto r_{\alpha,\beta}(h_{\nu})$ is strictly decreasing on $[0,\infty)$ and consequently we have $r_{\beta}^{c}(h_{\nu})<r_{\alpha,\beta}(h_{\nu})<r_{\beta }^{\ast}(h_{\nu})$ for all $\alpha\in(0,1),$ $\beta\in[0,1)$ and $\nu>-1.$
\end{Theorem}

\begin{figure}[!ht]
   \centering
       \includegraphics[width=10cm]{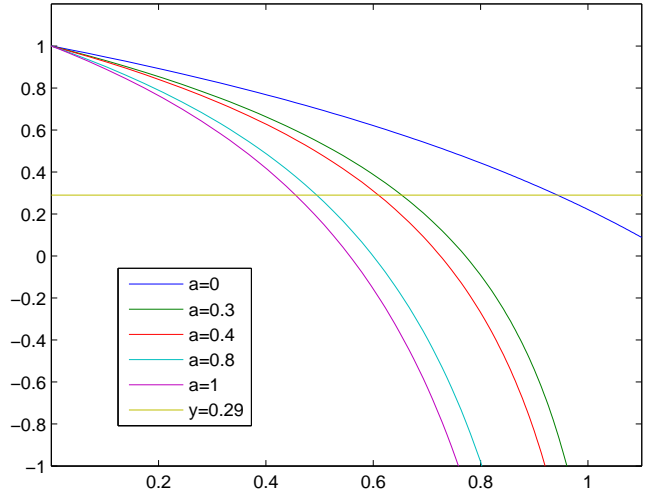}
       \caption{The graph of the function $r\mapsto \mathbb{J}(a,h_{-0.5}r))$ for $a\in\{0,0.3,0.4,0.8,1\}$ on $[0,1.1].$}
       \label{fig3}
\end{figure}

Finally , we mention that the cases $\alpha=0$ and $\alpha=1$ of the above Theorem were also considered recently in \cite[Theorem 1(c)]{bsk} and \cite[Theorem 1.3]{anal}. Our Theorem \ref{th3} is a common generalization of these results. Figure \ref{fig3} illustrates the fact that for $\alpha\in[0,1]$ the radius of $\alpha$-convexity of the function $h_{\nu}$ is between its radii of convexity and starlikeness, that is, $r_{\beta}^{c}(h_{\nu})<r_{\alpha,\beta}(h_{\nu})<r_{\beta }^{\ast}(h_{\nu})$ for all $\alpha\in(0,1),$ $\beta\in[0,1)$ and $\nu>-1.$ We considered the particular cases when $\beta=0.29,$ $\nu=-0.5$ and $\alpha\in\{0,0.3,0.4,0.8,1\}.$ We would like to take the opportunity to mention that \cite[Theorem 1(c)]{bsk} should be corrected as follows: if $\nu>-1,$ then $r^{\ast}(h_{\nu})=z_{\nu,\beta,1},$ where $z_{\nu,\beta,1}$ is the smallest positive root of the equation $z^{\frac{1}{2}}J_{\nu}'(z^{\frac{1}{2}})+(2-2\beta-\nu)J_{\nu}(z^{\frac{1}{2}})=0.$ In \cite[Theorem 1(c)]{bsk} the above result was stated wrongly with $z$ instead of $z^{\frac{1}{2}}.$ Consequently, \cite[Corollary 1(c)]{bsk} should be rewritten accordingly as follows: if  $\nu> -1,$ then the radius of starlikeness of the function $h_{\nu}$ is $z_{\nu,0,1},$ which
denotes the smallest positive root of the equation $z^{\frac{1}{2}}J_{\nu}'(z^{\frac{1}{2}})+(2-\nu)J_{\nu}(z^{\frac{1}{2}})=0.$

\section{\bf Proof of the main results}
\setcounter{equation}{0}

In this section our aim is to present the proofs of the main theorems. For convenience in the sequel we will use the following notation
$$\mathbb{J}(\alpha,u(z))=(1-\alpha)\frac{zu'(z)}{u(z)}+\alpha\left(1+\frac{zu''(z)}{u'(z)}\right).$$

\begin{proof}[\bf Proof of Theorem \ref{th1}]
Without loss of generality we assume that $\alpha>0.$ The case $\alpha=0$ was proved already in \cite{bsk}. By using the definition of the function $f_{\nu}$ we have
$$\frac{zf_{\nu}'(z)}{f_{\nu}(z)}=\frac{1}{\nu}\frac{zJ_{\nu}'(z)}{J_{\nu}(z)},\ \ \ \ 1+\frac{zf_{\nu}''(z)}{f'_{\nu}(z)}=1+\frac{zJ_\nu''(z)}{J_\nu'(z)}+\left(\frac{1}{\nu}-1\right)\frac{zJ_\nu'(z)}{J_\nu(z)}.$$
In view of the following infinite product representations \cite[p. 235]{nist}
$$J_{\nu}(z)=\frac{\left(\frac{1}{2}z\right)^{\nu}}{\Gamma(\nu+1)}\prod_{n\geq 1}\left(1-\frac{z^2}{j_{\nu,n}^2}\right), \ \
J_{\nu}'(z)=\frac{\left(\frac{1}{2}z\right)^{\nu-1}}{2\Gamma(\nu)}\prod_{n\geq 1}\left(1-\frac{z^2}{j_{\nu,n}'^2}\right),$$
where $j_{\nu,n}$ and $j_{\nu,n}'$ are the $n$th positive roots of $J_{\nu}$ and $J_{\nu}',$ respectively, logarithmic differentiation yields
$$\frac{zJ'_\nu(z)}{J_\nu(z)}=\nu-\sum_{n\geq1}\frac{2z^2}{j_{\nu,n}^2-z^2},\ \ 1+\frac{zJ_\nu''(z)}{J_\nu'(z)}=\nu-\sum_{n\geq1}\frac{2z^2}{j_{\nu,n}'^2-z^2},$$
which implies that
\begin{align*}
\mathbb{J}(\alpha,f_{\nu}(z))&=(1-\alpha)\frac{zf_{\nu}'(z)}{f_{\nu}(z)}+\alpha\left(1+\frac{zf_{\nu}''(z)}{f_{\nu}'(z)}\right)\\
&=\alpha+\left(\frac{1}{\nu}-\alpha\right)\frac{zJ_{\nu}'(z)}{zJ_{\nu}(z)}+\alpha\frac{zJ_{\nu}''(z)}{J_{\nu}'(z)}\\
&=1-\left(\frac{1}{\nu}-\alpha\right)\sum_{n\geq1}\frac{2z^2}{j_{\nu,n}^2-z^2}-\alpha\sum_{n\geq1}\frac{2z^2}{j_{\nu,n}'^2-z^2}.\end{align*}
On the other hand, we know (see \cite[Lemma 2.1]{anal}) that if $a>b>0,$ $z\in\mathbb{C}$ and $\lambda\leq 1,$ then for all $|z|<b$ we have
\begin{equation}\label{ineq}\lambda\real \left(\frac{z}{a-z}\right)-\real \left(\frac{z}{b-z}\right)\geq \lambda\frac{|z|}{a-|z|}-\frac{|z|}{b-|z|}.\end{equation} Note that in \cite[Lemma 2.1]{anal} it was assumed that $\lambda\in[0,1],$ however, following the proof of \cite[Lemma 2.1]{anal} it is clear that we do not need the assumption $\lambda\geq0.$ By using the inequality \eqref{ineq} for all $z\in\mathbb{D}(0,j_{\nu,1}')$ we obtain the inequality
\begin{equation*}\label{fnuconvex}\frac{1}{\alpha}\real\mathbb{J}(\alpha,f_{\nu}(z))\geq\frac{1}{\alpha}+\left(1-\frac{1}{\alpha\nu}\right)
\sum_{n\geq1}\frac{2r^2}{j_{\nu,n}^2-r^2}-\sum_{n\geq1}\frac{2r^2}{j_{\nu,n}'^2-r^2}=\frac{1}{\alpha}\mathbb{J}(\alpha,f_{\nu}(r)),\end{equation*}
where $|z|=r.$ Here we used that the zeros $j_{\nu,n}$ and $j_{\nu,n}'$ interlace according to the inequalities \cite[p. 235]{nist} \begin{equation}\label{interlace}\nu\leq j_{\nu,1}'<j_{\nu,1}<j_{\nu,2}'<j_{\nu,2}<j_{\nu,3}'<{\dots}.\end{equation} Now, the above deduced inequality implies that for $r\in(0,j_{\nu,1}')$ we have $\inf_{z\in\mathbb{D}(0,r)}\mathbb{J}(\alpha,f_{\nu}(z))=\mathbb{J}(\alpha,f_{\nu}(r)).$
On the other hand, the function $r\mapsto \mathbb{J}(\alpha,f_{\nu}(r))$ is strictly decreasing on $(0,j_{\nu,1}')$
since \begin{align*}\frac{\partial}{\partial r}\mathbb{J}(\alpha,f_{\nu}(r))&=-\left(\frac{1}{\nu}-\alpha\right)\sum_{n\geq1}\frac{4rj_{\nu,n}^2}{(j_{\nu,n}^2-r^2)^2}
-\alpha\sum_{n\geq1}\frac{4rj_{\nu,n}'^2}{(j_{\nu,n}'^2-r^2)^2}\\&<\alpha\sum_{n\geq1}\frac{4rj_{\nu,n}^2}{(j_{\nu,n}^2-r^2)^2}
-\alpha\sum_{n\geq1}\frac{4rj_{\nu,n}'^2}{(j_{\nu,n}'^2-r^2)^2}<0\end{align*}
for $\nu>0$ and $r\in(0,j_{\nu,1}').$ Here we used again that the zeros $j_{\nu,n}$ and $j_{\nu,n}'$ interlace and for all $n\in\mathbb{N},$ $\nu>0$ and $r<\sqrt{j_{\nu,1}j_{\nu,1}'}$ we have that $$j_{\nu,n}^2(j_{\nu,n}'^2-r^2)^2<j_{\nu,n}'^2(j_{\nu,n}^2-r^2)^2.$$ We also have that $\lim_{r\searrow0}\mathbb{J}(\alpha,f_{\nu}(r))=1>\beta$ and $\lim_{r\nearrow j_{\nu,1}'}\mathbb{J}(\alpha,f_{\nu}(r))=-\infty,$ which means that for $z\in\mathbb{D}(0,r_1)$ we have $\real\mathbb{J}(\alpha,f_{\nu}(z))>\beta$ if and only if $r_1$ is the unique root of $\mathbb{J}(\alpha,f_{\nu}(r))=\beta,$
situated in $(0,j_{\nu,1}').$ Finally, by using again the interlacing inequalities \eqref{interlace} we obtain the inequality
$$\frac{\partial}{\partial\alpha}\mathbb{J}(\alpha,f_{\nu}(r))=\sum_{n\geq1}\frac{2r^2}{j_{\nu,n}^2-r^2}-\sum_{n\geq1}\frac{2r^2}{j_{\nu,n}'^2-r^2}<0,$$
where $\nu>0,$ $\alpha\geq0$ and $r\in(0,j_{\nu,1}').$ This implies that the function $\alpha\mapsto \mathbb{J}(\alpha,f_{\nu}(r))$ is strictly decreasing on $[0,\infty)$ for all $\nu>0$ and $r\in(0,j_{\nu,1}')$ fixed. Consequently, as a function of $\alpha$ the unique root of the equation $\mathbb{J}(\alpha,f_{\nu}(r))=\beta$ is strictly decreasing, where $\beta\in[0,1),$ $\nu>0$ and $r\in(0,j_{\nu,1}')$ are fixed. Thus, in the case when $\alpha\in(0,1)$ the radius of $\alpha$-convexity of the function $f_{\nu}$ will be between the radius of convexity and the radius of starlikeness of the function $f_{\nu}.$ This completes the proof.
\end{proof}

\begin{proof}[\bf Proof of Theorem \ref{th2}]
Similarly, as in the proof of Theorem \ref{th1} we assume that $\alpha>0.$ The case $\alpha=0$ was proved already in \cite{bsk}. We start with the following relations
$$\frac{zg'_{\nu}(z)}{g_{\nu}(z)}=1-\nu+\frac{zJ'_{\nu}(z)}{J_{\nu}(z)},\ \ \ z\frac{g_{\nu}''(z)}{g_{\nu}'(z)}=\frac{\nu(\nu-1)J_\nu(z)+2(1-\nu)zJ_\nu'(z)+z^2J_\nu''(z)}{(1-\nu)J_\nu(z)+zJ_\nu'(z)}.$$
The recurrence formula \cite[p. 222]{nist} $zJ'_\nu(z)=\nu{J_\nu(z)}-zJ_{\nu+1}(z)$ and the fact that $J_{\nu}$ is a particular solution of the Bessel differential equation imply that
$$z\frac{g''_{\nu}(z)}{g'_{\nu}(z)}=z\frac{zJ_{\nu+2}(z)-3J_{\nu+1}(z)}{J_{\nu}(z)-zJ_{\nu+1}(z)},$$
and using \cite[Lemma 2.4]{anal} it follows that
$$1+z\frac{g''_{\nu}(z)}{g'_{\nu}(z)}=1-\sum_{n\geq1}\frac{2z^2}{\alpha_{\nu,n}^2-z^2},$$
where $\alpha_{\nu,n}$ is the $n$th positive zero of the Dini function $z\mapsto (1-\nu)J_\nu(z)+zJ'_\nu(z).$ Thus, we have that
\begin{align*}
\mathbb{J}(\alpha,g_{\nu}(z))&=(1-\alpha)\frac{zg_{\nu}'(z)}{g_{\nu}(z)}+\alpha\left(1+\frac{zg_{\nu}''(z)}{g_{\nu}'(z)}\right)\\
&=(1-\alpha)\left(1-\nu+\frac{zJ'_{\nu}(z)}{J_{\nu}(z)}\right)+\alpha\left(1+z\frac{zJ_{\nu+2}(z)-3J_{\nu+1}(z)}{J_{\nu}(z)-zJ_{\nu+1}(z)}\right)\\
&=1+(\alpha-1)\sum_{n\geq 1}\frac{2z^2}{j_{\nu,n}^2-z^2}-\alpha\sum_{n\geq 1}\frac{2z^2}{\alpha_{\nu,n}^2-z^2}.
\end{align*}
Applying the inequality \eqref{ineq} we have that
$$\frac{1}{\alpha}\real\mathbb{J}(\alpha,g_{\nu}(z))\geq\frac{1}{\alpha}+\left(1-\frac{1}{\alpha}\right)\sum_{n\geq 1}\frac{2r^2}{j_{\nu,n}^2-r^2}-\sum_{n\geq 1}\frac{2r^2}{\alpha_{\nu,n}^2-r^2}=\frac{1}{\alpha}\mathbb{J}(\alpha,g_{\nu}(r)),$$
where $|z|=r.$ Here we used tacitly that for all $n\in\{1,2,\dots\}$
we have $\alpha_{\nu,n}\in(j_{\nu,n-1},j_{\nu,n}),$ where $j_{\nu,
n}$ is the $n$th positive zero of $J_{\nu}.$ This follows immediately from Dixon's theorem \cite[p. 480]{watson}, which says that when $\nu>-1$ and
$a,b,c,d$ are constants such that $ad\neq bc,$ then the positive zeros of $z\mapsto aJ_{\nu}(z)+bzJ_{\nu}'(z)$ are interlaced with those of $z\mapsto cJ_{\nu}(z)+dzJ_{\nu}'(z).$ Thus, if we choose $a=1-\nu,$ $b=1,$ $c=1$ and $d=0,$ then the required assertion follows. Note also that the zeros $\alpha_{\nu,n}$ are all real when $\nu>-1,$ see \cite[p. 482]{watson}, and thus the application of the inequality \eqref{ineq} is allowed. Thus, for $r\in(0,\alpha_{\nu,1})$ we get $\inf_{z\in\mathbb{D}(0,r)}\real\mathbb{J}(\alpha,g_{\nu}(z))=\mathbb{J}(\alpha,g_{\nu}(r)),$ since according to the minimum principle of harmonic functions the infimum is taken on the boundary. On the other hand, the function $r\mapsto \mathbb{J}(\alpha,g_{\nu}(r))$ is strictly decreasing on $(0,\alpha_{\nu,1})$
since \begin{align*}\frac{\partial}{\partial r}\mathbb{J}(\alpha,g_{\nu}(r))&=(\alpha-1)\sum_{n\geq1}\frac{4rj_{\nu,n}^2}{(j_{\nu,n}^2-r^2)^2}
-\alpha\sum_{n\geq1}\frac{4r\alpha_{\nu,n}^2}{(\alpha_{\nu,n}^2-r^2)^2}\\&<\alpha\sum_{n\geq1}\frac{4rj_{\nu,n}^2}{(j_{\nu,n}^2-r^2)^2}
-\alpha\sum_{n\geq1}\frac{4r\alpha_{\nu,n}^2}{(\alpha_{\nu,n}^2-r^2)^2}<0\end{align*}
for $\nu>-1$ and $r\in(0,\alpha_{\nu,1}).$ Here we used again that the zeros $j_{\nu,n}$ and $\alpha_{\nu,n}$ interlace and for all $n\in\mathbb{N},$ $\nu>-1$ and $r<\sqrt{j_{\nu,1}\alpha_{\nu,1}}$ we have that $$j_{\nu,n}^2(\alpha_{\nu,n}^2-r^2)^2<\alpha_{\nu,n}^2(j_{\nu,n}^2-r^2)^2.$$ We also have that $\lim_{r\searrow0}\mathbb{J}(\alpha,g_{\nu}(r))=1>\beta$ and $\lim_{r\nearrow \alpha_{\nu,1}}\mathbb{J}(\alpha,g_{\nu}(r))=-\infty,$ which means that for $z\in\mathbb{D}(0,r_2)$ we have $\real\mathbb{J}(\alpha,g_{\nu}(z))>\beta$ if and only if $r_2$ is the unique root of $\mathbb{J}(\alpha,g_{\nu}(r))=\beta,$
situated in $(0,\alpha_{\nu,1}).$ Finally, by using again the interlacing inequalities $j_{\nu,n-1}<\alpha_{\nu,n}<j_{\nu,n}$ we obtain the inequality
$$\frac{\partial}{\partial\alpha}\mathbb{J}(\alpha,g_{\nu}(r))=\sum_{n\geq1}\frac{2r^2}{j_{\nu,n}^2-r^2}-\sum_{n\geq1}\frac{2r^2}{\alpha_{\nu,n}^2-r^2}<0,$$
where $\nu>-1,$ $\alpha\geq0$ and $r\in(0,\alpha_{\nu,1}).$ This implies that the function $\alpha\mapsto \mathbb{J}(\alpha,g_{\nu}(r))$ is strictly decreasing on $[0,\infty)$ for all $\nu>-1$ and $r\in(0,\alpha_{\nu,1})$ fixed. Consequently, as a function of $\alpha$ the unique root of the equation $\mathbb{J}(\alpha,g_{\nu}(r))=\beta$ is strictly decreasing, where $\beta\in[0,1),$ $\nu>-1$ and $r\in(0,\alpha_{\nu,1})$ are fixed. Thus, for $\alpha\in(0,1)$ the radius of $\alpha$-convexity of the function $g_{\nu}$ is between the radius of convexity and the radius of starlikeness of the function $g_{\nu}.$
\end{proof}

\begin{proof}[\bf Proof of Theorem \ref{th3}]
Similarly, as in the proof of Theorems \ref{th1} and \ref{th2} we assume that $\alpha>0.$ The case $\alpha=0$ was proved already in \cite{bsk}.
Combining
$$\frac{zh'_{\nu}(z)}{h_{\nu}(z)}=1-\frac{\nu}{2}+\frac{1}{2}\frac{z^{\frac{1}{2}}J'_{\nu}(z^{\frac{1}{2}})}{J_{\nu}(z^{\frac{1}{2}})}=1-\sum_{n\geq1}\frac{z}{j_{\nu,n}^{2}-z}$$
with \cite[Lemma 2.5]{anal}
$$z\frac{h_{\nu}''(z)}{h_{\nu}'(z)}=\frac{\nu(\nu-2)J_\nu(z^{\frac{1}{2}})+(3-2\nu)z^{\frac{1}{2}}J_\nu'(z^{\frac{1}{2}})
+zJ_\nu''(z^{\frac{1}{2}})}{2(2-\nu)J_\nu(z^{\frac{1}{2}})+2z^{\frac{1}{2}}J_\nu'(z^{\frac{1}{2}})}=-\sum_{n\geq1}\frac{z}{\beta_{\nu,n}^2-z},$$
where $\beta_{\nu,n}$ stands for the $n$th positive zero of the Dini function $z\mapsto (2-\nu)J_{\nu}(z)+zJ_{\nu}'(z),$ it follows that
\begin{align*}
\mathbb{J}(\alpha,h_{\nu}(z))=(1-\alpha)\frac{zh_{\nu}'(z)}{h_{\nu}(z)}+\alpha\left(1+\frac{zh_{\nu}''(z)}{h_{\nu}'(z)}\right)=1+(\alpha-1)\sum_{n\geq 1}\frac{z}{j_{\nu,n}^2-z}-\alpha\sum_{n\geq 1}\frac{z}{\beta_{\nu,n}^2-z}.
\end{align*}
Applying again the inequality \eqref{ineq} we have that
$$\frac{1}{\alpha}\real\mathbb{J}(\alpha,h_{\nu}(z))\geq\frac{1}{\alpha}+\left(1-\frac{1}{\alpha}\right)\sum_{n\geq 1}\frac{r}{j_{\nu,n}^2-r}-\sum_{n\geq 1}\frac{r}{\beta_{\nu,n}^2-r}=\frac{1}{\alpha}\mathbb{J}(\alpha,h_{\nu}(r)),$$
where $|z|=r.$ Here we used tacitly that for all $n\in\{1,2,\dots\}$
we have $\beta_{\nu,n}\in(j_{\nu,n-1},j_{\nu,n}),$ which follows immediately from Dixon's theorem \cite[p. 480]{watson}, similarly as in the case the roots $\alpha_{\nu,n}$ in the proof of Theorem \ref{th2}. Note also that the zeros $\beta_{\nu,n}$ are all real when $\nu>-1,$ see \cite[p. 482]{watson}, and thus the application of the inequality \eqref{ineq} is allowed. Thus, for $r\in(0,\beta_{\nu,1}^2)$ we get $\inf_{z\in\mathbb{D}(0,r)}\real\mathbb{J}(\alpha,h_{\nu}(z))=\mathbb{J}(\alpha,h_{\nu}(r)).$ On the other hand, the function $r\mapsto \mathbb{J}(\alpha,h_{\nu}(r))$ is strictly decreasing on $(0,\beta_{\nu,1}^2)$
since \begin{align*}\frac{\partial}{\partial r}\mathbb{J}(\alpha,h_{\nu}(r))&=(\alpha-1)\sum_{n\geq1}\frac{rj_{\nu,n}^2}{(j_{\nu,n}^2-r)^2}
-\alpha\sum_{n\geq1}\frac{r\beta_{\nu,n}^2}{(\beta_{\nu,n}^2-r)^2}\\&<\alpha\sum_{n\geq1}\frac{rj_{\nu,n}^2}{(j_{\nu,n}^2-r)^2}
-\alpha\sum_{n\geq1}\frac{r\beta_{\nu,n}^2}{(\beta_{\nu,n}^2-r)^2}<0\end{align*}
for $\nu>-1$ and $r\in(0,\beta_{\nu,1}^2).$ Here we used again that the zeros $j_{\nu,n}$ and $\beta_{\nu,n}$ interlace and for all $n\in\mathbb{N},$ $\nu>-1$ and $r<{j_{\nu,1}\beta_{\nu,1}}$ we have that $$j_{\nu,n}^2(\beta_{\nu,n}^2-r)^2<\beta_{\nu,n}^2(j_{\nu,n}^2-r)^2.$$ We also have that $\lim_{r\searrow0}\mathbb{J}(\alpha,h_{\nu}(r))=1>\beta$ and $\lim_{r\nearrow \beta_{\nu,1}}\mathbb{J}(\alpha,g_{\nu}(r))=-\infty,$ which means that for $z\in\mathbb{D}(0,r_3)$ we have $\real\mathbb{J}(\alpha,h_{\nu}(z))>\beta$ if and only if $r_3$ is the unique root of $\mathbb{J}(\alpha,h_{\nu}(r))=\beta,$
situated in $(0,\beta_{\nu,1}^2).$ Finally, by using again the interlacing inequalities $j_{\nu,n-1}<\beta_{\nu,n}<j_{\nu,n}$ we obtain the inequality
$$\frac{\partial}{\partial\alpha}\mathbb{J}(\alpha,h_{\nu}(r))=\sum_{n\geq1}\frac{r}{j_{\nu,n}^2-r}-\sum_{n\geq1}\frac{r}{\beta_{\nu,n}^2-r}<0,$$
where $\nu>-1,$ $\alpha\geq0$ and $r\in(0,\beta_{\nu,1}^2).$ This implies that the function $\alpha\mapsto \mathbb{J}(\alpha,h_{\nu}(r))$ is strictly decreasing on $[0,\infty)$ for all $\nu>-1$ and $r\in(0,\beta_{\nu,1}^2)$ fixed. Consequently, as a function of $\alpha$ the unique root of the equation $\mathbb{J}(\alpha,h_{\nu}(r))=\beta$ is strictly decreasing, where $\beta\in[0,1),$ $\nu>-1$ and $r\in(0,\beta_{\nu,1}^2)$ are fixed. Thus, when $\alpha\in(0,1)$ the radius of $\alpha$-convexity of the function $h_{\nu}$ is between the radius of convexity and the radius of starlikeness of the function $h_{\nu}.$
\end{proof}

\end{document}